\documentclass[12pt]{article}

\usepackage{tikz,color,url}
\usepackage{amsmath,amssymb,latexsym}

\newtheorem{theorem}{Theorem}[section]

\newtheorem{lemma}[theorem]{Lemma}
\newtheorem{proposition}[theorem]{Proposition}

\newcommand{\proof}{\noindent{\bf Proof.\ }}
\newcommand{\qed}{\hfill $\square$ \medskip}

\newcommand{\dimnl}{{\rm dim}_{\rm n\ell}}
\newcommand{\diml}{{\rm dim_{\ell}}}
\newcommand{\diam}{{\rm diam}}

\newcommand{\ex}{{\rm ex}}

\begin{document}

\title{Nonlocal metric dimension of graphs}

\author{
Sandi Klav\v zar $^{a,b,c}$
\and
Dorota Kuziak $^d$
}

\date{}

\maketitle

\begin{center}
$^a$ Faculty of Mathematics and Physics, University of Ljubljana, Slovenia \\
\medskip

$^b$ Faculty of Natural Sciences and Mathematics, University of Maribor, Slovenia \\
\medskip

$^{c}$ Institute of Mathematics, Physics and Mechanics, Ljubljana, Slovenia \\
{\tt sandi.klavzar@fmf.uni-lj.si}
\medskip

$^{d}$ Departamento de Estad\'istica e Investigaci\'on Operativa, Universidad de C\'adiz, Algeciras, Spain \\
{\tt dorota.kuziak@uca.es}
\end{center}

\begin{abstract}
Nonlocal metric dimension ${\rm dim}_{\rm n\ell}(G)$ of a graph $G$ is introduced as the cardinality of a smallest nonlocal resolving set, that is, a set of vertices which resolves each pair of non-adjacent vertices of $G$. Graphs $G$ with ${\rm dim}_{\rm n\ell}(G) = 1$ or with ${\rm dim}_{\rm n\ell}(G) = n(G)-2$ are characterized. The nonlocal metric dimension is determined for block graphs, for corona products, and for wheels. Two upper bounds on the nonlocal metric dimension are proved. An embedding of an arbitrary graph into a supergraph with a small nonlocal metric dimension and small diameter is presented. 
\end{abstract}


\noindent
{\bf AMS Subj.\ Class.\  (2020)}: 05C12, 05C76

\section{Introduction}
\label{sec:intro}

Let $G = (V(G), E(G))$ be a graph, $X\subseteq V(G)$, and $u,v\in V(G)$. Then $u$ and $v$ are {\em resolved} by $X$ if there exists $x\in X$ such that $d_G(u,x) \ne d_G(v,x)$, where $d_G(y,z)$ denotes the shortest-path distance between vertices $y$ and $z$ of $G$. A set $X$ such that each pair of vertices of $G$ is resolved by $X$ is a {\em resolving set} of $G$. A smallest resolving set is a {\em metric basis} of $G$, its cardinality is the {\em metric dimension} of $G$ denoted by $\dim(G)$. This concept is one of the central ones in the field of metric graph theory, and its popularity stems, among other things, from the fact that it has found many applications in fields as diverse as computer science, chemistry, biology, and social sciences. For a better understanding of the concept and its applications see the recent survey~\cite{till-2022+}, while for a comprehensive survey on its variants see the other recent survey \cite{dorota-2022+}.

In 2010, Okamoto, Phinezy, and Zhang~\cite{okamoto-2010} introduced the local metric dimension of a graph as a natural version of the metric dimension. In this version, which was further researched in particular in~\cite{abrishami-2022, barr-2019, barr-2016, fancy-2017, klavzar-2020,  rodriguez-2016, rodriguez-2015, susi-2015}, we need to resolve only pairs of adjacent vertices. The dual concept, in which all pairs of non-adjacent vertices are to be resolved, has surprisingly not been considered in the literature. In this paper we fill this gap.

If $X\subseteq V(G)$ resolves each pair of non-adjacent vertices, then we speak of a {\em nonlocal resolving set}. A smallest nonlocal resolving set is a {\em nonlocal metric basis} of $G$, its cardinality is the {\em nonlocal metric dimension} of $G$ and will be denoted by $\dimnl(G)$. By definition,
\begin{equation}
\label{eq:dim-dimnl}
\dimnl(G) \le \dim(G)
\end{equation}
holds for every graph $G$. The difference can be arbitrary large as already demonstrated by complete graphs for which we have $\dimnl(K_n) = 0$ and $\dim(K_n) = n-1$ for every $n\ge 1$.

We proceed as follows. In the next section we first recall definitions, notation, and a result needed later. Then we characterize the graphs with the nonlocal metric dimension equal to $1$, and show that $\dimnl(G) = \dim(G)$ holds for bipartite graphs $G$. In Section~\ref{sec:block-corona} we determine the nonlocal metric dimension for block graphs and for corona products. The latter result is reduced to a join of a graph with $K_1$. This  motivates us to determine the nonlocal metric dimension of wheels in Section~\ref{sec:wheels}. In the subsequent section we prove two upper bounds on the nonlocal metric dimension and characterize graphs $G$ with $\dimnl(G) = n(G) - 2$. In Section~\ref{sec:embed}, we embed an arbitrary graph into a supergraph with a small nonlocal metric dimension and small diameter. In the concluding section directions for further study are indicated.  

\section{Preliminaries}
\label{sec:prelim}

Unless stated otherwise, graphs considered will be connected. The order of a graph $G$ will be denoted by $n(G)$. If $X=\{x_1,\dots,x_k\}$ and $u\in V(G)$, then the {\em metric representation of u with respect to $X$} is the vector $r(u|X)=(d_G(u,x_1),\dots, d_G(u,x_k))$. The diameter $\diam(G)$ of $G$ is the largest distance between pairs of vertices of $G$. The clique number of $G$ is the order of a largest complete subgraph in $G$ and denoted by $\omega(G)$. As usual, the chromatic number of $G$ is denoted by $\chi(G)$. The edge cover number $\beta'(G)$ of $G$ is  the smallest number of edges such that each vertex is incident with at least one of these edges. The join $G+H$ of graphs $G$ and $H$ is obtained from the disjoint union of a copy of $G$ and a copy of $H$ by adding an edge between each vertex of $G$ and each vertex of $H$. The complement of $G$ will be denoted by $\overline{G}$. For a positive integer $k$ we will use the notation $[k] = \{1,\ldots, k\}$.

A vertex of degree at least~$3$ in a tree $T$ is called a \emph{branch vertex}. A leaf $u$ of $T$ is called a \emph{terminal leaf} of a branch vertex $v$ of $T$ if $d_T(u,v) < d_T(u,w)$ for every other branch vertex $w$ of $T$.  A branch vertex $v$ of $T$ is an \emph{exterior branch vertex} of $T$ if it has at least one terminal leaf. The path from an exterior branch vertex to its terminal leaf is called a \emph{terminal path}. Let $n_1(T)$ denote the number of leaves of $T$, and let $\ex(T)$ denote the number of exterior branch vertices of $T$. In~\cite{chartrand-2000, HaMe-76, Slater-75} it was proved that if $T$ is a tree that is not a path, then
\begin{equation}
\label{eq:dim-for-trees}
\dim(T) = n_1(T) - \ex(T).
\end{equation}

We already mentioned that $\dimnl(K_n) = 0$ for $n\ge 1$. Moreover, as soon as $G$ is not complete, $\dimnl(G)\ge 1$. Clearly, $\dimnl(P_n) = 1$ for $n\ge 3$. It is also not difficult to see that $\dimnl(C_n) = 2$ for $n\ge 4$. 

To characterize graphs $G$ with $\dimnl(G) = 1$, we need the following concept. If $x\in V(G)$ and $k\ge 0$, then $L_k(x) = \{u\in V(G):\ d_G(x,u) = k\}$ is a {\em distance level} of $x$. In particular, $L_0(x) = \{x\}$, and $L_1(x)$ is the (open) neighborhood of $x$.

\begin{proposition}
\label{prop:dim=1}
If $G$ is a non-complete graph, then $\dimnl(G) = 1$ if and only if there exists a vertex $x$ such that $L_k(x)$ induces a complete graph for every $k \le \diam(G)$.
\end{proposition}

\proof
Assume first that $\dimnl(G) = 1$ and let $\{x\}$ be a nonlocal metric basis. If $u, v\in L_k(x)$ for some $k\ge 1$, then $d_G(x,u) = d_G(x,v) = k$ and hence $u$ and $v$ must be adjacent. Thus $L_k(x)$ induces a complete graph. ($L_0(x)$ is trivially complete.)

Conversely, let $x$ be a vertex whose all distance levels induce complete graphs. Therefore, if $u, v\in V(G)$ and $uv\notin E(G)$, then $d_G(x,u) \ne d_G(x,v)$ and so $u$ and $v$ are resolved by $x$.
\qed

In general $\dimnl(G)$ can be arbitrary smaller than $\dim(G)$. This cannot happen if $G$ is bipartite as the next result asserts. 

\begin{proposition}
\label{prop:dimnl-bipartite}
If $G$ is a bipartite graph with $n(G)\ge 3$, then $\dimnl(G) = \dim(G)$.
\end{proposition}

\proof
By~\eqref{eq:dim-dimnl} we have $\dimnl(G) \le \dim(G)$. Let $X$ be a nonlocal metric basis of $G$. As $n(G)\ge 3$, $G$ is not complete and hence $|X| \ge 1$. 
Consider arbitrary vertices $u,v \in V(G)$. If $uv \notin E(G)$, then $u$ and $v$ are resolved by $X$ because $X$ is a nonlocal resolving set. If $uv \in E(G)$, then $u$ and $v$ are resolved by $X$ since $G$ is bipartite. It follows that $X$ is also a resolving set. Hence $\dim(G) \le |X| = \dimnl(G)$ and we are done.
\qed

\section{Block graphs and corona products}
\label{sec:block-corona}

Proposition~\ref{prop:dimnl-bipartite} implies that if $T$ is a tree with $n(T)\ge 3$, then $\dimnl(T) = \dim(T)$. This fact generalizes to block graphs as follows. Let $G$ be a block graph. The {\em block-cutpoint tree} $\widehat{G}$ of $G$ is the tree whose vertices are the blocks and the cut-vertices of $G$, and a block $B$ is adjacent to a cut-vertex $v$ if $v\in V(B)$, cf.~\cite[Definition~4.1.20]{west-2001}. 

\begin{theorem}
\label{thm:block-graphs}
If $G$ is a block graph with $n(G)\ge 3$, then $\dimnl(G) = \dim(\widehat{G})$.
\end{theorem}

\proof
We can extend the terminology for trees from the previous section to block graphs as follows. Let $G$ be a block graph. A block $B$ of $G$ is called a \emph{branch block} if the vertices of $B$ have at least three independent neighbors outside of $B$. A branch block $B$ of $G$ is an \emph{exterior branch block} of $G$ if in $\widehat{G}$ there exist $t\ge 1$ terminal paths emanating from a cut-vertex of $B$. If so, we say that $B$ has $t$ {\em rays}.  

Let now $X$ be a nonlocal resolving set of $G$. If $B$ is an exterior branch block of $G$ with $t$ rays, then  we claim that $X$ has at least $t-1$ vertices in these rays, each one from a different ray. Indeed, otherwise there exist two rays $R_1$ and $R_2$ with no vertex from $X$. Let $R_1$ and $R_2$ be attached to $B$ at vertices $u_1$ and $u_2$ (note that $u_1 = u_2$ is possible), and select neighbors $v_1$ and $v_2$ of $u_1$ and $u_2$ in $R_1$ and $R_2$, respectively, Then $v_1$ and $v_2$ are not adjacent and have the same metric representation with respect to $X$.  This proves the claim.
Each exterior branch block of $G$ corresponds to an exterior branch vertex of $\widehat{G}$. Since $\widehat{G}$ is a tree, \eqref{eq:dim-for-trees} implies that $\dimnl(G) \ge \dim(\widehat{G})$. 

To prove that $\dimnl(G) \le \dim(\widehat{G})$ holds, let $X$ be a metric basis of $\widehat{G}$. We may without loss of generality assume that every element of $X$ is a leaf of $\widehat{G}$. Then each element $x\in X$ corresponds to a terminal block $B_x$ of $G$, that is, a block which contains exactly one cut-vertex. For every $x\in X$, let $w_x$ be an arbitrary but fixed vertex of $B_x$ which is not a cut-vertex. Set $Y = \{w_x:\ x\in X\}$. Since we do not need to resolve adjacent vertices, it is straightforward to check that $Y$ is a nonlocal resolving set of $G$. Hence $\dimnl(G) \le |Y| = |X| = \dim(\widehat{G})$. 
\qed

The metric dimension of block graphs was studied in~\cite{hoffmann-2016}. It was proved that  the metric dimension of a block graph $G$ is equal to the metric dimension of a tree, obtained from $G$ by replacing each block $B$ of size at least $3$ by a star whose leaves are the vertices of $B$. The local metric dimension of block graphs was investigated in~\cite{rodriguez-2015}. 

Let $G$ and $H$ be graphs where $V(G) = \{g_1, \ldots,g_{n(G)}\}$.  The \emph{corona product} of $G$ and $H$, denoted by $G\odot H$, is a graph obtained from the disjoint union of a copy of $G$ and $n(G)$ copies of $H$, denoted by $H_i$, $i\in [n(G)]$. The product $G\odot H$ is then constructed by making $g_i$ adjacent to every vertex in $H_i$ for each $i\in [n(G)]$.  Let further $\widetilde{H_i}$ be the subgraph of $G\odot H$ induced by $V(H_i)\cup \{g_i\}$. Clearly, $\widetilde{H_i}\cong H + K_1$. We will use this notation in the rest of the section.

 The metric dimension of corona products has been investigated in~\cite{iswadi-2011, yero-2011}. The local metric dimension of corona products has been studied in~\cite{rodriguez-2016} and more generally of generalized hierarchical products in~\cite{klavzar-2020}. In particular, in~\cite{rodriguez-2016} it is proved for the local metric dimension $\diml$ that if $G$ is a connected graph and $H$ is a graph of radius at least $4$, then $\diml(G\odot H)=n(G)\cdot \diml(H+K_1)$. In general, $\diml(G\odot H)\leq n(G)\cdot \diml(H+K_1)$, see~\cite{klavzar-2020}. Here we add the following formula for the nonlocal metric dimension.

\begin{theorem}
\label{thm:corona}
If $G$ is a graph and $H$ a non-complete graph, then 
$$\dimnl(G\odot H) = n(G)\cdot \dimnl(H + K_1).$$
\end{theorem}

\proof
Let $X$ be a nonlocal metric basis of $G\odot H$ and let $X_i = X \cap V(\widetilde{H_i})$. If $u,v\in V(H_i)$ and $w\in V(G\odot H)\setminus V(H_i)$, then $d_{G\odot H}(u,w) = d_{G\odot H}(v,w)$. Since also $d_{G\odot H}(u,g_i) = d_{G\odot H}(v,g_i) = 1$, any two non-adjacent vertices from $H_i$ must be resolved by some vertex from $H_i$. Therefore, $X_i$ is a nonlocal resolving set of $\widetilde{H_i}$ and so $|X_i| \ge \dimnl(\widetilde{H_i}) = \dimnl(H + K_1)$. Consequently,
$$\dimnl(G\odot H) = |X| = \sum_{i=1}^{n(G)} |X_i| \ge \sum_{i=1}^{n(G)} \dimnl(\widetilde{H_i}) = n(G)\cdot \dimnl(H + K_1).$$

Let now $Y$ be a nonlocal metric basis of $H\odot K_1$. Since $H$ is not complete, $|Y|\ge 1$. For each $i\in [n(G)]$, let $Y_i$ be the copy of $Y$ in $\widetilde{H_i}$. Note that $g_i\notin Y_i$. We claim that $\cup_{i=1}^{n(G)}Y_i$ is a nonlocal resolving set of $G\odot H$. Indeed, if vertices $u$ and $v$ from  $V(H_i)$ are non-adjacent, then they are resolved by some vertex from $Y_i$.  Assume next $u\in V(\widetilde{H_i})$ and $v\in V(\widetilde{H_j})$,  where $i\ne j$. If $y_i \in Y_i$, then $d(y_i,u) \leq 2$ and $d(y_i,v) \geq 3$, hence $u$ and $v$ are again resolved. We have thus seen that $\cup_{i=1}^{n(G)}Y_i$ is a nonlocal resolving set and henceforth,
$$\dimnl(G\odot H) \le |\cup_{i=1}^{n(G)}Y_i| = n(G)\cdot |Y| = n(G)\cdot \dimnl(H + K_1).$$
We conclude that $\dimnl(G\odot H) = n(G)\cdot \dimnl(H + K_1)$.
\qed

In the case when the second factor of a corona product is complete, we can bound the nonlocal metric dimension as follows. 

\begin{theorem}
\label{thm:corona-complete}
If $G$ is a graph and $n\ge 1$, then 
$$\dim(G) \le \dimnl(G\odot K_n) \le n(G).$$
\end{theorem}

\proof
To prove the upper bound we claim that $V(G)$ is a nonlocal resolving set of $G\odot K_n$. Indeed, let $u, v \in V(G\odot K_n)\setminus V(G)$ be two non-adjacent vertices. Then $u\in H_i$ and $v\in H_j$ for some $i$ and $j$, where $i\ne j$. As $d_{G\odot K_n}(u,g_i) = 1$ and $d_{G\odot K_n}(v,g_i)\ge 2$, the claim follows. Hence $\dimnl(G\odot K_n) \le n(G)$. 

Let now $X$ be a nonlocal metric basis of $G\odot H$; then $|X| = \dimnl(G\odot H)$. Let $Y\subseteq V(G)$ be the set obtained from $X$ by replacing each vertex $u\in X\cap V(\widetilde{H_i})$ by $g_i$. Note that this in particular means that if $g_i\in X$, then also $g_i\in Y$. Clearly, $|Y| \le |X|$. We claim that $Y$ is a resolving set of $G$. For this sake consider arbitrary vertices $g_i, g_j\in V(G)\subseteq V(G\odot H)$. Let $u\in V(H_i)$ and $v\in V(H_j)$. Since $X$ is a nonlocal metric basis of $G\odot H$ and $uv\notin E(G\odot H)$, there exists a vertex $w\in X$ such that $d_{G\odot H}(u,w) \ne d_{G\odot H}(v,w)$. If $w\in \widetilde{H_i} \cup \widetilde{H_j}$, then $g_i\in Y$ or $g_j\in Y$ and there is nothing to be proved. Suppose henceforth that $w\in \widetilde{H_k}$, where $k\ne i,j$. If $w\ne g_k$, then
$$1 + d_G(g_i, g_k) + 1 = d_{G\odot H}(u,w) \ne d_{G\odot H}(v,w) = 1 + d_G(g_j, g_k) +1,$$
and if $w = g_k$, then
\begin{align*}
1 + d_G(g_i, g_k) & = 1 + d_G(g_i, w) \\
& = d_{G\odot H}(u,w) \ne d_{G\odot H}(v,w) \\
& = 1 + d_G(g_j, w) = 1 + d_G(g_j, g_k).
\end{align*} 
In either case we get $d_G(g_i, g_k) \ne d_G(g_j, g_k)$. As $g_k\in Y$, we have proved that $Y$ is a resolving set of $G$ which in turn implies that
$$\dim(G) \le |Y| \le |X| = \dimnl(G\odot H)$$
and we are done.
\qed

If $r,s\ge 3$, then it can be checked that $\dimnl(K_{r,s}\odot K_n) = r + s - 2 = \dimnl(K_{r,s}) = \dim(K_{r,s})$. This shows the sharpness of the lower bound of Theorem~\ref{thm:corona-complete}. On the other hand, if $m\ge 2$ and $n\ge 1$, then $\dimnl(P_m\odot K_n) = 2$. Since $\dim(P_m) = 1$ and $n(K_m) =m$, the case $m=2$ shows the tightness of the upper bound, while the cases $m\ge 3$ demonstrate that the intermediate values are also attainable. 

\section{Nonlocal metric dimension of wheels}
\label{sec:wheels}

In Theorem~\ref{thm:corona}, the nonlocal metric dimension of a corona product $G\odot H$ is reduced to the nonlocal metric dimension of the join of $H$ with $K_1$. To determine the latter deserves a special attention. Here we solve it for the join of a cycle with $K_1$, that is, for wheels. Recall that the wheel graph $K_1+C_{n}$ of order $n+1\ge 4$ is denoted by $W_{1,n}$.

As a consequence of~\cite[Corollary 5(iv)]{rodriguez-2016} we obtain $\diml(W_{1,n}) = \lceil n/4\rceil$. The metric dimension of $W_{1,n}$ was independently studied in~\cite{buczkowski-2003, shan-2002}, where it is proved that if $n\ge 7$, then
$$\dim(W_{1,n}) = \left\lfloor{\frac{2n+2}{5}}\right\rfloor.$$
Clearly, $\dimnl(W_{1,3})=1$ and $\dimnl(W_{1,4})=\dimnl(W_{1,5})=\dimnl(W_{1,6})=2$. The main result of this section reads as follows.

\begin{theorem}
\label{thm:dimnl-wheel}
If $n\ge 7$, then
$$\dimnl(W_{1,n})=\left\lfloor{\frac{2n}{5}}\right\rfloor.$$
\end{theorem}

Note that if $n\ge 7$, then $\dim(W_{1,n}) = \dimnl(W_{1,n})$ if and only if $n \bmod 5 \in \{0, 1, 3\}$. Otherwise, $\dimnl(W_{1,n}) = \dim(W_{1,n}) - 1$. In the rest of the section we prove Theorem~\ref{thm:dimnl-wheel}, for which some additional terminology is needed.

Let $V(C_n)=\{0,1,\dots,n-1\}$ with natural adjacency. Operation with vertices will be done  modulo $n$. Let $v$ be the central vertex of the wheel $W_{1,n}$, where $n\ge 7$. Notice first that, $v$ does not belong to any nonlocal metric basis, since for any $0\le i\le n-1$, $v$ and $i$ are adjacent. Let $X\subset V(C_n)$ be a set of vertices such that $|X|\ge 2$. A {\em gap} of $X$ is a set of vertices $A_{i,j}=\{i+1,\dots,j-1\}$, $|i-j|\ge 1$, of $C_n$ such that $i,j\in X$ and $\{i+1,\dots,j-1\}\cap X=\emptyset$. We call $i$ and $j$ {\em neighboring vertices} of $X$. Two gaps $A_{i,j}$ and $A_{j,k}$ are called {\em neighboring gaps}.

\begin{lemma}
\label{lem:dimnl-wheel-gaps}
Let $n\ge 7$, and let $S$ be a nonlocal metric basis of $W_{1,n}$. Then the following hold.
\begin{enumerate}
  \item[(i)] $|A_{i,j}| \le 4$ for each gap of $S$.
  \item[(ii)] $|A_{i,j}| \ge 3$ holds for at most one gap of $S$.
  \item[(iii)] If $|A_{i,j}|\ge 2$, then each of the two neighboring gaps of $A_{i,j}$ contains at most one vertex.
  \item[(iv)] Moreover, $W_{1,n}$ contains a nonlocal metric basis $S'$ that fulfils conditions (i)-(iii) and in addition has no gap of size $3$.
\end{enumerate}
\end{lemma}

\proof
(i) Consider the gap $A_{i,j}$ and suppose that $j\ge i+6$, where $i\le n-1$. Then the vertices $i+2$ and $i+4$ have the same metric representation with respect to $S$. As they are not adjacent, we have a contradiction.

(ii) Suppose there exist two distinct gaps $A_{i,i+k}$ and $A_{j,j+k'}$, where $k, k'\in \{4, 5\}$. Then for two non-adjacent vertices $i+2$ and $j+2$ we have $r(i+2|S)=r(j+2|S)$, a contradiction.

(iii) Let  $A_{i,j}$ be a gap with $j\ge i+3$. Consider now a neighboring gap, without loss of generality let it be $A_{i',i}$, and suppose that $i - i' > 2$. Then the vertices $i-1$ and $i+1$ have the same metric representation with respect to $S$. Again, as they are not adjacent, we have a contraction.

(iv) Assume that $S$ contains a gap $A_{i,i+4}$. By (ii), this gap is the only gap of size $3$. Moreover, by (iii), there are the following cases to be considered. Assume first that $i+5\notin S$. Then $i+6$ must belong to $S$. Let $S' = (S \cup \{i+5\})\setminus \{i+4\}$. Having in mind that $i+2$ and $i+3$ are adjacent, it is straightforward to verify that $S'$ is a required nonlocal metric basis. The case when $i-1\notin S$ is treated analogously. The last case to consider is when $i-1, i+5\in S$. Then $i+6\notin S$ and the set $S' = (S \cup \{i+6\})\setminus \{i+4\}$ is a required nonlocal metric basis.
\qed

\begin{lemma}
\label{lem:dimnl-wheel-set}
Let $n\ge 7$. If $X\subseteq V(C_n)$ satisfies conditions (i)-(iii) of Lemma~\ref{lem:dimnl-wheel-gaps}, then $X$ is a nonlocal resolving set of $W_{1,n}$.
\end{lemma}

\proof
Notice that by Proposition~\ref{prop:dim=1}, $\dimnl(W_{1,n})\ge 2$, for $n\ge 7$. Let $u\in V(W_{1,n})\setminus X$. As the central vertex $v$ is adjacent to any other vertex, we may assume that $u\ne v$. We then distinguish three cases.

\medskip\noindent
\textbf{Case 1}: $u$ belongs to a gap $A_{i,i+2}$ (gap of size $1$). \\
Then $u=i+1$ and $u$ has distance $1$ to $i,i+2\in X$ and distance $2$ to all other vertices of $X$. Hence any other vertex has different metric representation with respect to $X$ because $n\ge 7$.

\medskip\noindent
\textbf{Case 2}: $u$ belongs to a gap $A_{i,i+3}$ (gap of size $2$). \\
Let without loss of generality $u=i+1$. Then $u$ has distance $1$ to $i$ and distance $2$ to all other vertices of $X$. By condition (iii) from Lemma~\ref{lem:dimnl-wheel-gaps}, there is no any other vertex with the same metric representation.

\medskip\noindent
\textbf{Case 3}: $u$ belongs to a gap $A_{i,i+j}$, $j\in \{4,5\}$. \\
If $u=i+1$, then $u$ has distance $1$ to $i$ and distance $2$ to all the other vertices of $X$. Hence, by condition (iii) of Lemma~\ref{lem:dimnl-wheel-gaps}, there is no other vertex with the same metric representation. The case when $u=i+j-1$ is done similarly. Now, if $u\in \{i+2,i+j-2\}$, then $r(u|X)=(2,\dots,2)$. If $j=4$, then there is only one such vertex, and if $j=5$, there are two such vertices. In the latter case these two vertices are adjacent. By condition (ii) from Lemma~\ref{lem:dimnl-wheel-gaps}, no other vertex has this metric representation.
\qed

We are now in a position to prove Theorem~\ref{thm:dimnl-wheel}. We construct a set $S\subseteq V(W_{1,n})$ with $|S|=\left\lfloor{\frac{2n}{5}}\right\rfloor$ depending on  the residue class modulo $5$.

\medskip\noindent
{\bf Case 1}: $n \bmod 5\in \{0,1\}$. \\
Then $n=5k$ or $n=5k+1$, where $k\ge 2$. In both cases let $S=\{5i,5i+2: i\in [k-1]\}\cup \{0,5k-1\}$. Notice that $|S|=2k=\left\lfloor{\frac{2n}{5}}\right\rfloor$.

\medskip\noindent
{\bf Case 2}: $n\equiv 2 \bmod 5$. \\
Then $n=5k+2$, where $k\ge 1$. Let $S=\{5i,5i+2: 1\le i\le k-1\}\cup \{0,5k\}$. Then $|S|=2k=\left\lfloor{\frac{2n}{5}}\right\rfloor$.

\medskip\noindent
{\bf Case 3}: $n \bmod 5 \in \{3,4\}$. \\
Then $n=5k+3$ or $n=5k+4$, where $k\ge 1$. In both cases let $S=\{5i,5i+2: i\in [k]\}\cup \{0\}$. Then $|S|=2k+1=\left\lfloor{\frac{2n}{5}}\right\rfloor$.

\medskip
In each case, the set $S$ fulfils the conditions of Lemma~\ref{lem:dimnl-wheel-set}, hence $S$ is a nonlocal resolving set of $W_{1,n}$ and thus $\dimnl(W_{1,n})\le\left\lfloor{\frac{2n}{5}}\right\rfloor$.

To prove that $\dimnl(W_{1,n}) \ge \left\lfloor{\frac{2n}{5}}\right\rfloor$, consider an arbitrary nonlocal resolving set $X$ of $W_{1,n}$. In view of Lemma~\ref{lem:dimnl-wheel-gaps}(iv) we may without loss of generality assume that $X$ has no gaps of size $3$. We distinguish two cases.

\medskip\noindent
{\bf Case 1}: $|X|=2r$, where $r\ge 1$. \\
By Lemma~\ref{lem:dimnl-wheel-gaps}(iii), at most $r$ gaps contain two or four vertices, and by Lemma~\ref{lem:dimnl-wheel-gaps}(ii), at most one of these gaps contains $4$ vertices. Hence, there is at most $2(r-1) + 4 + r = 3r+2$ vertices belonging to the gaps. Since $n\le 2r+(3r+2)$, we have that $|X|=2r \ge 2\left(\frac{n-2}{5}\right)$. Since $2r$ is an integer we conclude that $|X| \ge \left\lceil{\frac{2}{5}n-\frac{4}{5}}\right\rceil = \left\lfloor{\frac{2n}{5}}\right\rfloor$.

\medskip\noindent
{\bf Case 2}: $|X|=2r+1$, where $r\ge 1$. \\
Again, By Lemma~\ref{lem:dimnl-wheel-gaps}(iii), at most $r$ gaps contain two or four vertices, and by Lemma~\ref{lem:dimnl-wheel-gaps}(ii), at most one of these gaps contains $4$ vertices. Hence in this case there is at most $2(r-1) + 4 + (r+1) = 3r+3$ vertices belonging to the gaps. Since $n\le 2r+1+(3r+3)$, we have that $|X|=2r+1\ge \left\lceil{\frac{2}{5}n-\frac{3}{5}}\right\rceil \ge \left\lfloor{\frac{2n}{5}}\right\rfloor$.

\medskip
This proves Theorem~\ref{thm:dimnl-wheel}.

\section{Upper bounds}

In this section we prove two upper bounds on the nonlocal metric dimension. Along the way we show that complete bipartite graphs can be characterized as the graphs $G$ with $\dimnl(G) = n(G) - 2$. 

\begin{proposition}
\label{prop:omega}
If $G$ is a graph, then $\dimnl(G) \le n(G) - \omega(G)$. In particular, if $n(G)\ge 2$, then $\dimnl(G) \le n(G) - 2$.
\end{proposition}

\proof
Since $\dimnl(K_n) = 0$ for $n\ge 1$, the result holds for complete graphs. Hence assume in the rest that $G$ is not complete and $n(G) \ge 3$. 

Let $Q$ be a clique of $G$ of order $\omega(G)$. Then we claim that $V(G) \setminus V(Q)$ is a nonlocal resolving set. Indeed, let $u$ and $v$ be non-adjacent vertices. At least one of them is not in $Q$, say $u$. Then $d_G(u,u) = 0 < d_G(u,v)$ and since $u\in V(G) \setminus V(Q)$, the claim is proved. Hence $\dimnl(G) \le n(G) - n(Q) = n(G) - \omega(G)$. The second assertion, $\dimnl(G) \le n(G) - 2$, follows using the same argument except that we consider an arbitrary $K_2$ as a complete subgraph. 
\qed

In view of Proposition~\ref{prop:omega} we next describe graphs that attain the bound $n(G) - 2$.

\begin{proposition}
\label{prop:n(G)-2}
If $G$ is a graph, then $\dimnl(G) = n(G) - 2$ if and only if $G = K_{s,t}$, where $s\ge 1$ and $t\ge 2$.
\end{proposition}

\proof
If $\dimnl(G) = n(G) - 2$, then by~\eqref{eq:dim-dimnl} we have $\dim(G)\ge n(G)-2$. Since the only graphs $G$ with $\dim(G) = n(G)-1$ are complete graphs, we can restrict our attention to the graphs  $G$ with  $\dim(G) = n(G)-2$. In~\cite{chartrand-2000} it is proved that $\dim(G) = n(G) - 2$ if and only if $G$ is one of the following graphs: $K_{s,t}$ ($s\ge 1, t\ge 2$), $K_s+ \overline{K_t}$ ($s\ge 1, t\ge 2$), and $K_s+ (K_1\cup K_t)$ ($s\ge 1, t\ge 1$). By~\eqref{eq:dim-dimnl} we know that in each of these cases we have $\dimnl(G) \le n(G) - 2$. Hence we need to determine in which of these cases we have $\dimnl(G) \ge n(G) - 2$.

Consider first complete bipartite graphs $K_{s,t}$ and let $S$ and $T$ be its bipartition sets, where $|S| = s$ and $|T| = t$. If $u,v\in S$, then $u$ and $v$ have the same distance to all the other vertices. If follows that for each nonlocal metric basis $X$ we have $|X\cap S| = s-1$. Analogously we get $|X\cap T| = t-1$. Hence $|X| = (s-1) + (t-1) = n(K_{s,t}) - 2$.

Consider next joins $K_s+ \overline{K_t}$, $s\ge 1$, $t\ge 2$. The case $s=1$ has already been treated in the above paragraph, so assume that $s, t\ge 2$. Then $\omega(K_s+ \overline{K_t}) = s + 1 \ge 3$, and hence by Proposition~\ref{prop:omega} we have $\dimnl(K_s+ \overline{K_t}) \le n(K_s+ \overline{K_t}) - \omega(K_s+ \overline{K_t}) < n(K_s+ \overline{K_t}) - 2$.

The last class of graphs to be considered is $K_s+ (K_1\cup K_t)$, $s\ge 1$, $t\ge 1$. If $s=t=1$, then $G=P_3$ which has already been considered earlier. And if $s\ge 2$ or $t\ge 2$, then $\omega(G)\ge 3$ and we conclude  as in the previous case that none of the graphs qualifies for the theorem.
\qed

\begin{theorem}
\label{thm:girth5}
If $G$ is a graph of girth at least $7$, then $\dimnl(G) \le \beta'(G) - 1$. Moreover, if $G$ is a tree, then the equality holds if and only if $G$ is obtained from a star by subdividing all but one of the edges at most once.
\end{theorem}

\proof
Note first that the assertion is true for $K_2$ since $\dimnl(K_2) = 0$. It is also clear that the result holds for paths since $\dimnl(P_n) = 1$ for $n\ge 3$. Hence in the rest  $n(G)\ge 4$ and $\Delta(G)\ge 3$.

Assume first that there exists an edge cover $S$ of $G$ with $|S| = \beta'(G)\ge 2$, such that $S$ contains edges $xy$ and $xy'$, where $y\ne y'$. We distinguish two cases.

\medskip\noindent
{\bf Case 1}: $\deg_G(y) = \deg_G(y') = 1$. \\
Then $\deg_G(x) \ge 3$. We define a set $X\subset V(G)$ as follows. Put into $X$ all the neighbors of $x$ but $y'$. Complete the construction of $X$ by putting into it an arbitrary vertex from each of the edges $f\in S$ which has no vertex yet in $X$. Let $x'$ be a neighbor of $x$ different from $y$ and $y'$. Since $G$ is triangle-free, the vertex $x'$ is covered either by the edge $xx'$, or by an edge $x'x''$, where $x''$ is not adjacent to $x$. Therefore, $|X| \le \beta'(G) - 1$. Note also that $x\notin X$. We claim that $X$ is a nonlocal resolving set. For this sake select arbitrary vertices $u,v\notin X$ such that $uv\notin E(G)$.  If $u\in \{x, y'\}$, then because $\deg_G(y) = 1$, the vertices $u$ and $v$ are resolved by $y$. Assume next that $u, v\notin \{x, y'\}$. Then there exist vertices $u'$ and $v'$ such that $uu', vv'\in E(G)$ and $u',v'\in X$. We claim that $u$ and $v$ are resolved by $u'$ and $v'$. If this is not the case, then since $d_G(u,u') = 1$ we have $d_G(v,u') = 1$, and because $d_G(v,v') = 1$ we have $d_G(u,v') = 1$. But this implies that $uu'vv'u$ is a $4$-cycle, a contradiction with our girth assumption.

\medskip\noindent
{\bf Case 2}: $\deg_G(y) \ge 2$. \\
Let $y''$ be a neighbor of $y$ different from $x$. Then we infer that $yy''\notin S$. Indeed, if $yy''\in S$, then $S' = S \setminus \{xy\}$ is an edge cover of $G$ of cardinality $\beta'(G) - 1$. Hence $y''$ is an end-vertex of some edge $y''y'''\in S$. Define now a set $X\subset V(G)$ as follows. First put $x$ and $y''$  into $X$. In addition, for any edge $f\in S\setminus \{xy, xy', y''y'''\}$ put into $X$ an arbitrary vertex from $f$. Note that $|X| \le \beta'(G) - 1$. We claim that $X$ is a nonlocal resolving set. Let $u,v\notin X$ such that $uv\notin E(G)$. Then we see that in every case there exist vertices $u'$ and $v'$ such that $uu', vv'\in E(G)$ and $u',v'\in X$. In particular, if $u=y$ and $v=y'$, then set $u'=y''$ and $v'=x$. Similarly, if $u=y$ and $v=y'''$, then select $u' = x$ and $v' = y''$. And if $u=y'$ and $v=y'''$, then set $u' = x$ and $v' = y''$. In any case we see that the vertices $u$, $v$, $u'$ and $v'$ lead to a $4$-cycle.

\medskip
Assume second that every edge cover $G$ has cardinality $n(G)/2$. Let $S$ be such an edge cover. Again we distinguish two cases.

\medskip\noindent
{\bf Case 1}: every edge from $S$ has an end-vertex of degree $1$. \\
Let $S = \{x_iy_i:\ i\in [n(G)/2]\}$ and assume without loss of generality that $\deg(x_i) = 1$ for all $i\in [n(G)/2]$. Set $X = \{x_2, \ldots, x_{n(G)/2}\}$. We claim that $X$ is a nonlocal resolving set. Let $u$ and $v$ be arbitrary non-adjacent vertices from $V(G)\setminus X$. Then at least one of $u$ or $v$ is from $\{y_2, \ldots, y_{n(G)/2}\}$, say $u=y_i$, $i\ge 2$. But then $x_i$ resolves $u$ and $v$, and the claim is proved. 

\medskip\noindent
{\bf Case 2}: there is an edge $xy\in S$ such that $\deg_G(x)\ge 2$ and $\deg_G(y)\ge 2$.  \\
Let $x'$ be a neighbor of $x$ different from $y$, and let $y'$ be a neighbor of $y$ different from $x$. By the assumption on $S$, there exist vertices $x''$ and $y''$ such that $x'x'', y'y''\in S$. Let $X\subset V(G)$ be the set containing $x'$, $y'$, and an arbitrary vertex from each edge $f\in S' = S\setminus \{xy, x'x'', y'y''\}$. We claim that $X$ is a nonlocal resolving set. Let $f=ww'\in S'$, where $w'\in X$. Then since $d_G(w,w') = 1$ and since each of the vertices $x$, $y$, $x''$, and $y''$ has either $x'\in X$ or $y'\in X$ as a neighbor, we get that $X$ resolves $w$ from each of the vertices from $\{x,y,x'',y''\}$. Similarly we infer that the pairs $x'',y$, and $x'',y''$, and $x,y''$ are resolved by $X$. The only pairs left to be considered are $x,x''$ and $y,y''$. Suppose on the contrary that the pair $x,x''$ is not resolved by $X$. Since $d_G(x,y') = 2$, we must have $d_G(x'',y') = 2$. Then $y'x'\in E(G)$, or $y''x''\in E(G)$, or there is a new vertex $z$ such that $y'z, zx''\in E(G)$. In either case we have a cycle of length at most $6$ in $G$ which is not possible. The other pair $y,y''$ is treated analogously.

\medskip
Let now $G$ be a tree. Let $W$ be the set of exterior branch vertices of $T$, and for a vertex $w\in W$, let $p(w)$ be the number of terminal paths of $w$. Then we have
$$\beta'(T) \ge \sum_{w\in W} p(w).$$
On the other hand, by~\eqref{eq:dim-for-trees} and Proposition~\ref{prop:dimnl-bipartite} we have
$$\dimnl(T) = \sum_{w\in W} (p(w)-1).$$
Hence $\dimnl(G) = \beta'(G) - 1$ can hold only when $|W| = 1$. Suppose hence that $G$ is a tree with only one exterior branch vertex. Then $G$ is obtained from a star $K_{1,n}$, $n\ge 3$, by subdividing each of its edges an arbitrary number of times. By~\eqref{eq:dim-for-trees}, $\dimnl(G) = n-1$. Moreover, if at least one edge of $K_{1,n}$ has been subdivided at least two times, then $\beta'(G) \ge n+1$. Suppose next that each edge of $K_{1,n}$ has been subdivided exactly once. Then $\beta'(G) = n+1$. So we are left with trees $G$  obtained from $K_{1,n}$ by subdividing all but one of the edges at most once. In this case we have $\dimnl(G) = n-1 = \beta'(G) - 1$.
\qed

\section{Embeddings}
\label{sec:embed}

In our final result we prove that any graph can be embedded as an induced subgraph into a supergraph with a small nonlocal metric dimension and small diameter. More precisely, the following holds. 

\begin{theorem}
\label{thm:embed}
Every connected graph $G$ is an induced subgraph of a graph $H$ with $\dimnl(H) \le \left\lceil\log(\chi(\overline{G}))\right\rceil$ and $\diam(H) \le 4$. Moreover, if $2^{k-1} \le \chi(\overline{G})< 2^k$, for some integer $k$, then $\diam(H) \le 3$.
\end{theorem}

\proof
If $\chi(\overline{G}) = 1$, then $\overline{G}$ is an edge-less graph and hence $G$ is complete. Then, by definition, $\dimnl(G) = 0$ and hence the conclusion holds in this case by embedding the complete graph into itself. We may thus assume in the rest that $\chi(\overline{G}) = s \ge 2$. Let $k$ be the unique integer with $2^{k-1} < s \le 2^{k}$.

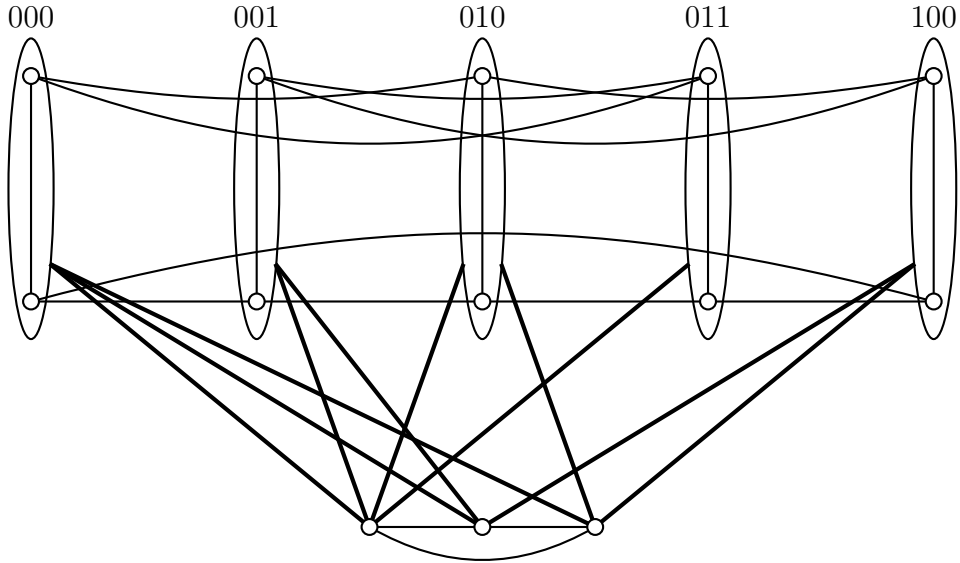
\begin{figure}[ht!]
\begin{center}
\begin{tikzpicture}[scale=1.0,style=thick]
\tikzstyle{every node}=[draw=none,fill=none]
\def\vr{3pt} 

\begin{scope}[yshift = 0cm, xshift = 0cm]
\path (0,0) coordinate (x1);
\path (3,0) coordinate (x2);
\path (6,0) coordinate (x3);
\path (9,0) coordinate (x4);
\path (12,0) coordinate (x5);
\path (0,3) coordinate (y1);
\path (3,3) coordinate (y2);
\path (6,3) coordinate (y3);
\path (9,3) coordinate (y4);
\path (12,3) coordinate (y5);
\path (4.5,-3) coordinate (z1);
\path (6,-3) coordinate (z2);
\path (7.5,-3) coordinate (z3);
\draw (x1) -- (x2) -- (x3) -- (x4) -- (x5);
\draw (x1) to[bend left=15] node {} (x5);
\draw (y1) to[bend left=-10] node {} (y3);
\draw (y3) to[bend left=-10] node {} (y5);
\draw (y2) to[bend left=-10] node {} (y4);
\draw (y1) to[bend left=-20] node {} (y4);
\draw (y2) to[bend left=-20] node {} (y5);
\draw (0,1.5) ellipse (0.3cm and 2.0cm);
\draw (3,1.5) ellipse (0.3cm and 2.0cm);
\draw (6,1.5) ellipse (0.3cm and 2.0cm);
\draw (9,1.5) ellipse (0.3cm and 2.0cm);
\draw (12,1.5) ellipse (0.3cm and 2.0cm);
\foreach \i in {1,...,5}
\draw (x\i) -- (y\i);
\draw (z1) -- (z2) -- (z3);
\draw (z1) to[bend left=-30] node {} (z3);

\draw[ultra thick] (z1) to (0.25,0.5);
\draw[ultra thick] (z2) to (0.25,0.5);
\draw[ultra thick] (z3) to (0.25,0.5);
\draw[ultra thick] (z1) to (3.25,0.5);
\draw[ultra thick] (z2) to (3.25,0.5);
\draw[ultra thick] (z1) to (5.75,0.5);
\draw[ultra thick] (z3) to (6.25,0.5);
\draw[ultra thick] (z1) to (8.75,0.5);
\draw[ultra thick] (z2) to (11.75,0.5);
\draw[ultra thick] (z3) to (11.75,0.5);
\draw (x1)  [fill=white] circle (\vr);
\draw (x2)  [fill=white] circle (\vr);
\draw (x3)  [fill=white] circle (\vr);
\draw (x4)  [fill=white] circle (\vr);
\draw (x5)  [fill=white] circle (\vr);
\draw (y1)  [fill=white] circle (\vr);
\draw (y2)  [fill=white] circle (\vr);
\draw (y3)  [fill=white] circle (\vr);
\draw (y4)  [fill=white] circle (\vr);
\draw (y5)  [fill=white] circle (\vr);
\draw (z1)  [fill=white] circle (\vr);
\draw (z2)  [fill=white] circle (\vr);
\draw (z3)  [fill=white] circle (\vr);
\draw[above] (y1)++(0.0,0.5) node {$000$};
\draw[above] (y2)++(0.0,0.5) node {$001$};
\draw[above] (y3)++(0.0,0.5) node {$010$};
\draw[above] (y4)++(0.0,0.5) node {$011$};
\draw[above] (y5)++(0.0,0.5) node {$100$};

\end{scope}

\end{tikzpicture}
\end{center}
\caption{Embedding the Petersen graph into a graph $H$ as in the proof of Theorem~\ref{thm:embed}. A thick edge means that a bottom vertex is adjacent to both encircled vertices.}
\label{fig:graph-H}
\end{figure}

Let $X_0, \ldots, X_{s-1}$ be the color classes of $\overline{G}$ under some optimal coloring of it. By our assumption, $2^{k-1} < s \le 2^{k}$. Then the sets $X_0, \ldots, X_{s-1}$ respectively induce complete subgraphs of $G$. Since $G$ is connected, we may without loss of generality assume that there exists an edge connecting a vertex from $X_0$ by a vertex from $X_{s-1}$. Construct now a graph $H$ as follows. First take the disjoint union of $G$ and $K_k$. For each $X_i$ write $i$ in its binary representation, say $i = i_1\ldots i_k$. Then for each $j$ such that $i_j = 0$, add all the edges between the vertex $j\in V(K_k)$ and $X_i$. See Fig.~\ref{fig:graph-H} for an example of the construction where $G$ is the Petersen graph $P$. For it note that $\chi(\overline{P}) = 5$. 

We claim that $D = V(K_k)$ is a nonlocal resolving set of $H$. For this sake let $u$ and $v$ be arbitrary non-adjacent vertices of $H$. If at least one of $u$ and $v$ belongs to $D$, there is nothing to prove. Otherwise, $u \in X_i$ and $v\in X_j$, where $i\ne j$. If $i = i_1\ldots i_k$ and $j = j_1\ldots j_k$ are the binary representations of $i$ and $j$, then there exists $p$ such that, without loss of generality, $i_p = 0$ and $j_p = 1$. Then $d_{H}(p,u) = 1 \ne 2 = d_{H}(p,v)$, hence $u$ and $v$ are resolved. This proves the claim. Since $|D| = k = \left\lceil\log(s)\right\rceil = \left\lceil\log(\chi(\overline{G}))\right\rceil$, we conclude that $\dimnl(H) \le \left\lceil\log(\chi(\overline{G}))\right\rceil$.

To show that $\diam(H) \le 4$, let again $u$ and $v$ be arbitrary vertices of $H$. Note first that if $u\in V(K_k)$, then $d_{H}(u,v)\le 2$. Assume hence that $u\in X_i$ and $v\in X_j$, where $i < j$. Suppose first $j < 2^k$. Then each of $u$ and $v$ has a neighbor in $V(K_k)$ which implies that $d_{H}(u,v)\le 3$. Note that this fact in particular implies that if $\chi(\overline{G}) < 2^{k}$, then $\diam(H) \le 3$. Suppose second that $j = 2^k$. By our selection of $X_0$ and $X_{2^k}$, there is an edge $v'w$, where $v'\in X_{2^k}$ and $w\in X_0$. As $w$ is adjacent to all the vertices of $V(K_k)$ and $u$ has at least one neighbor in $V(K_k)$, we have a $u,v$-path of length $4$.
\qed

A result for the standard metric dimension which would be parallel to Theorem~\ref{thm:embed} is not possible. The reason is the following. Imagine a connected, non-complete graph $G$ and having a color class (in an optimal coloring of it) of order $2^r$, $r\ge 2$. Then $G$ contains a clique of the same order and hence in whichever graph $H$ we embed $G$ as an induced subgraph, we have $\omega(H)\ge 2^r$. Then~\cite[Theorem 1]{chartrand-2000} which asserts that if $G$ is a graph with $\diam(G) = d$, then $f(n,d) \le \dim(G) \le n - d$, where $f(n,d)$ is the least positive integer $k$ for which $k + d^k \ge n$, implies that $\dim(H)\ge r$. So $\dim(H)$ can be arbitrary larger than $\left\lceil\log(\chi(\overline{G}))\right\rceil$ (actually, arbitrary larger than $\chi(\overline{G})$ for that matter).

\section{Concluding remarks}

In this article we have introduced the concept of nonlocal metric dimension, which seems to deserve wider interest, not least because it is in some sense a complementary concept to the established local metric dimension. Many problems about this new concept remain open, let us mention here just a few that follow from our results.

In Theorem~\ref{thm:corona} we have proved a formula for the local metric dimension of corona products with the second factor being not complete. Extending this theorem (as done in~\cite{klavzar-2020} for the local metric dimension) to generalized hierarchical products is open. In view of Theorem~\ref{thm:corona-complete} it would be interesting to prove a formula for $\dimnl(G\odot K_n)$. Theorem~\ref{thm:girth5} is proved for graphs of girth at least $7$. We suspect that the condition can be relaxed to girth at least $5$. Moreover, the edge cover number in triangle-free graphs (that is, of girth at least $4$) is equal to the clique cover number. Hence extending Theorem~\ref{thm:girth5} to arbitrary graphs, where the upper bound would be stated as a function of the clique covering number would be of interest.

\section*{Acknowledgements}

We would like to express our sincere thanks to the reviewers for their very careful reading of the paper.

Sandi Klav\v{z}ar acknowledges the financial support from the Slovenian Research Agency through research core funding No.\ P1-0297 and projects J1-2452 and N1-0285. Dorota Kuziak was partially supported by the Spanish Ministry of Science and Innovation through the grant PID2019-105824GB-I00. and by ``Plan Propio de Investigaci\'{o}n-UCA" ref.\ no.\ EST2022-075.

\end{document}